\documentclass[conference]{IEEEtran}

\usepackage{pifont}
\usepackage{bbding}  

\IEEEoverridecommandlockouts                              

\usepackage{amsmath,amssymb,epsf,epsfig,times}
\usepackage{amsfonts}
\usepackage[all]{xy}
\usepackage{color}
\usepackage{subfigure}
\usepackage{url,cite}

\newtheorem{theorem}{Theorem}[section]
\newtheorem{lemma}{Lemma}[section]

\def\QED{\mbox{\rule[0pt]{1.5ex}{1.5ex}}}
\def\endproof{\hspace*{\fill}~\QED\par\endtrivlist\unskip}
\newcommand{\re}{\mathbb{R}}

\newtheorem{assumption}[theorem]{Assumption}

\newtheorem{remark}[theorem]{Remark}

\newcommand{\OMIT}[1]{}

\newif\ifpdf
\ifx\pdfoutput\undefined \pdffalse \else \pdfoutput=1 \pdftrue \fi
\ifpdf
\usepackage{graphicx}
\usepackage{epstopdf}
\usepackage{graphicx}
\fi

\title{\LARGE \bf
	Distributed Average Tracking for Second-order Agents with Nonlinear Dynamics}

\author{Sheida Ghapani, Salar Rahili, and Wei Ren 
	\thanks{The authors are with Department of Electrical Engineering, University of California, Riverside, CA
		92521, USA.
		{Email: sghap001@ucr.edu, srahi001@ucr.edu, ren@ee.ucr.edu}.
		This work was supported by National Science Foundation under Grant CMMI-1537729.}
}

\begin{document}

\maketitle
\thispagestyle{empty}
\pagestyle{empty}

\begin{abstract}		
 This paper addresses distributed average tracking of physical second-order agents with nonlinear dynamics, where the interaction among the agents is described by an undirected graph.
In both agents' and reference inputs' dynamics, there is a nonlinear term that satisfying the Lipschitz-type condition. 
To achieve the distributed average tracking problem in the presence of nonlinear term, a non-smooth filter and a control input are designed for each agent. 
The idea is that each filter outputs converge to the average of the reference inputs and the reference velocities asymptotically and in parallel each agent's position and velocity are driven to track its filter outputs. 
To overcome the nonlinear term unboundedness effect, novel state-dependent time varying gains are employed in each agent's filter and control input.
In the proposed algorithm, each agent needs its neighbors' filters outputs besides its own filter outputs, absolute position and absolute velocity and its neighbors' reference inputs and reference velocities.
Finally, the algorithm is simplified to achieve the distributed average tracking of physical second-order agents in the presence of an unknown bounded term in both agents' and reference inputs' dynamics. 
	\end{abstract}

	\section{INTRODUCTION}
	
	In this paper,  the distributed average tracking problem for a team of agents is studied, where each agent uses local information to calculate the average of individual time varying reference inputs, one per agent.
	The problem has found applications in distributed sensor fusion \cite{olfati2004consensus}, feature-based map merging \cite{aragues2012distributed}, distributed optimization \cite{salarsheida2015}, and distributed Kalman filtering \cite{bai2011}, where the scheme has been mainly used as an estimator.
	However, there are some applications such as region following formation control \cite{cheah2009region} and coordinated path planning \cite{vsvestka1998coordinated} that require the agents' physical
	states instead of estimator
	states to converge to a time varying network quantity, where each agent only has a local and incomplete copy of that quantity.
	Since the desired trajectory is not available to any agent, the distributed average tracking poses more theoretical challenges than the consensus and tracking problem.
	
	In the literature, some researchers have employed linear distributed algorithms, where the time varying reference inputs are required to satisfy restrictive constraints
	\cite{spanos2005dynamic,freeman2006,bai2010,kiaauthority,kiaDACsingularity}. 
	In \cite{freeman2006}, a proportional algorithm and a proportional-integral algorithm are proposed to achieve distributed average tracking for slowly-varying
	reference inputs with a bounded tracking
	error, where accurate estimator initialization is relaxed in the proportional-integral algorithm.
	In \cite{bai2010}, an extension of the proportional integral algorithm is employed for a special group of time varying reference inputs with a common denominator in their Laplace transforms, where the denominator is required to be used in the estimator  design.
	In \cite{kiaauthority}, a distributed average tracking problem is addressed with some steady-state errors, where the privacy of each agent's input is preserved.

	However, linear algorithms cannot ensure distributed average tracking for a more general group of reference inputs. Therefore, some researchers employ nonlinear tracking algorithms.
	In \cite{Nosrati2012}, a class of nonlinear algorithms is proposed for reference inputs with bounded deviations, where the tracking error is proved to be bounded.
	A non-smooth algorithm is proposed in \cite{chen2012distributed}, which is able to track arbitrary time varying reference inputs with bounded derivatives.
	However, all the aforementioned works study the distributed average tracking problem from a distributed estimation perspective without the requirement for agents to obey certain physical dynamics.
	There are various applications, where the distributed average tracking problem is relevant for designing distributed control laws for physical agents.
	The region-following formation control is one application \cite{cheah2009region}, where a swarm of robots move inside a dynamic region while keeping a desired formation.
	Distributed average tracking for double-integrator agents is studied in \cite{feidoubleintegrator}, where the reference inputs are allowed to have bounded accelerations.
	Refs.  \cite{zhao2013distributed} and \cite{FeiRobust} study the distributed average tracking for physical agents with general linear dynamics, where reference inputs and their control inputs are bounded.
In particular,  \cite{FeiRobust} proposes a discontinuous algorithm, while a continuous algorithm is employed in \cite{zhao2013distributed} with, respectively, static and adaptive coupling strengths.
		Ref. \cite{ghapani2015distributed} introduces a discontinuous algorithm and filter for a group of physical double-integrator agents, where each agent uses the relative positions and neighbors' filter outputs to remove the velocity measurements.
		
	However, in real applications physical agents might need to track the average of a group of real time varying reference inputs, where both physical agents and reference inputs have more complicated dynamics rather than single- or double-integrator dynamics.
	In fact, the dynamics of the physical agents and reference inputs must be taken into account in the control law design and the dynamics themselves introduce further challenges in the tracking and stability analyses.
		Therefore, the control law designed for physical agents with single- and double-integrator dynamics can no longer be used directly for physical agents subject to more complicated dynamic equations.	
	For example, \cite{fei2015EulerDAC} extends a proportional-integral control scheme  to achieve distributed average tracking for physical Euler-Lagrange systems for two different kinds of reference inputs with steady states and with bounded derivatives.  

		
In this paper, a distributed algorithm (controller design combined with filter design) is introduced to achieve the distributed average tracking for physical second-order agents with nonlinear dynamics.
Here, a nonlinear term satisfying the Lipshitz-type condition is considered in both agents' and reference inputs' dynamics to describe more complicated dynamics.
Due to the presence of the nonlinear term in the agents' dynamics, a local filter is introduced for each agent to estimate the average of the reference inputs and the reference velocities.
Further, a non-smooth control input is introduced to drive the agents to track the filter outputs.
Since the unknown term can be unbounded, we are faced with extra challenges. 
Therefore, novel time varying state-dependent gains are proposed in each agent's filter and control input to overcome the  unboundedness effect.
Finally, the filter dynamics and control input are simplified with constant gains to achieve the distributed average tracking, where the nonlinear term is bounded.

{\it Notations:} Throughout the paper, $\mathbb{R}$ denotes the set of all real numbers and $\mathbb{R}^+$ the set of all positive real numbers. 
	The transpose of matrix $A$ and vector $x$ are shown as $A^T$ and $x^T$, respectively. 
	Let $\mathbf{1}_n$ and $\mathbf{0}_n$ denote, respectively, the $n \times 1$ column vector of all ones and all zeros.
	Let $\mbox{diag}(z_1,\ldots,z_p)$ be the diagonal matrix with diagonal entries $z_1$ to $z_p$.
	We use 
	$\otimes$ to denote the Kronecker product, and $\mbox{sgn}(\cdot)$ to denote the $\mbox{signum}$ function defined component-wise. 
	For a vector function ${x(t):\re\mapsto\re^m}$, define $\|x\|_\mathfrak{p}$ as the $\mathfrak{p}$-norm, 
	${x(t)\in\mathbb{L}_2}$ if
	$\int_{0}^{\infty} \|x(\tau)\|^2_2\mbox{d}\tau<\infty$ and
	${x(t)\in\mathbb{L}_{\infty}}$ if for each element of $x(t)$, denoted as
	$x_i(t)$, ${\sup_{t \geq 0}|x_i(t)|<\infty}$, $i=1,\ldots,m$.
\section{Problem Statement}
Consider a multi-agent system consisting of $n$ physical agents described by the following nonlinear second-order dynamics	
\begin{align} \label{non-agents-dynamic}
\dot{x}_i(t)=&v_i(t), \notag  \\
\dot{v}_i(t)=&f(x_i(t),v_i(t),t)+u_i(t),  \qquad i=1,\ldots,n,
\end{align}
where $x_i(t) \in \mathbb{R}^\mathfrak{p}$, $v_i(t) \in \mathbb{R}^\mathfrak{p}$ and $u_i(t) \in \mathbb{R}^\mathfrak{p}$ are $i$th agent's position, velocity and control input, respectively, and  $f:\mathbb{R}^\mathfrak{p} \times \mathbb{R}^\mathfrak{p} \times \mathbb{R}^+ \to \mathbb{R}^\mathfrak{p}$ is a 
	vector-valued term which will be defined later.
	
	
	An \textit{undirected} graph $G \triangleq (V,E)$ is used to characterize the interaction topology among the agents, where ${V \triangleq \{1,\ldots,n\}}$ is the node set and $E \subseteq V \times V$ is the edge set.
	An edge $(j,i) \in E$ means that node $i$ can obtain information from node $j$ and vice versa.
	Self edges $(i,i)$ are not considered here.
	Let $m$ denote the number of edges in $E$, where the edges $(j,i)$ and $(i,j)$ are counted only once.
	The \textit{adjacency matrix} ${\mathbf{A} \triangleq [a_{ij}] \in \mathbb{R}^{n \times n}}$ of the graph $G$ is defined such that the edge weight ${a_{ij}=1}$ if ${(j,i) \in E}$ and ${a_{ij}=0}$ otherwise. For an undirected graph, ${a_{ij}=a_{ji}}$.
	The \textit{Laplacian matrix} ${L \triangleq [l_{ij}] \in \mathbb{R}^{n \times n}}$ associated with $\mathbf{A}$ is defined as ${l_{ii}=\sum_{j \ne i} a_{ij}}$ and ${l_{ij}=-a_{ij}}$, where ${i \ne j}$.
	For an undirected graph, $L$ is symmetric positive semi-definite.
	By arbitrarily assigning an orientation for the edges in $G$, let $D \triangleq [d_{ij}] \in  \mathbb{R}^{n \times m}$ be the \textit{incidence matrix} associated with $G$, where $d_{ij} = -1$ if the edge $e_j$ leaves node $i$, $d_{ij} = 1$ if it enters node $i$, and $d_{ij} = 0$ otherwise.
	The \textit{Laplacian matrix} $L$ is then given by $L=DD^T$ \cite{GodsilRoyle01}.
	\begin{assumption} \label{conn-graph}
		Graph $G$ is connected.
	\end{assumption}
	\begin{lemma} \cite{GodsilRoyle01} \label{eigen}
		Under Assumption \ref{conn-graph}, the \textit{Laplacian matrix} $L$ has a simple zero eigenvalue such that $0=\lambda_1(L)<\lambda_2(L) \leq \ldots \leq \lambda_n(L)$, where $\lambda_i(\cdot)$ denotes the $i$th eigenvalue. Furthermore, for any vector $y \in \mathbb{R}^n$ satisfying ${\mathbf{1}_n^T y=0}$, $\lambda_2(L) y^Ty \leq y^T L y \leq \lambda_n(L) y^Ty$.
	\end{lemma}
	Suppose that each agent has a time varying reference input $r_i(t) \in \mathbb{R}^\mathfrak{p}$, $i=1,\ldots,n$, satisfying
	\begin{align*} 
	\dot{r}_i(t)=& v_i^r(t),   \notag \\
	\dot{v}_i^r(t)=& f(r_i,v_i^r,t),
	\end{align*}
	where $v_i^r(t) \in \mathbb{R}^\mathfrak{p}$ is the reference velocity.
	Define ${r(t)=[r^T_1,\ldots,r^T_n ]^T}$,  ${v^r(t)=[{v_1^r}^T,\ldots,{v_n^r}^T ]^T}$ and ${f(r,v^r,t)=[f^T(r_1,v^r_1,t),\ldots,f^T(r_n,v^r_n,t)]^T}$.
	
Here the goal is to design $u_i(t)$ for agent $i$, $i=1,\ldots,n$, to track the average of the reference inputs and reference velocities, i.e.,
	\begin{align}\label{goal}
	\lim \limits_{t \to \infty} \|x_i(t)-\frac{1}{n} \sum_{j=1}^n r_j(t)\|_2=&0, \notag \\
	\lim \limits_{t \to \infty} \|v_i(t)-\frac{1}{n} \sum_{j=1}^n v_j^r(t)\|_2=&0, 
	\end{align}
where each agent has only local interaction with its neighbors. 	
As it was mentioned, there are many applications that the physical agents should track a time varying trajectory, where each agent has an incomplete copy of this trajectory.
While, the physical agents and reference trajectory might be described by more complicated dynamics rather than the double-integrator dynamics in real applications.
Therefore, we consider a more general group of physical agents, where the unknown term $f(\cdot,\cdot,t)$ in their dynamics satisfies the Lipschitz-type condition. 
\subsection{Main Result
	}\label{section2}
In this subsection, we study the distributed average tracking of second-order multi-agent systems with nonlinear dynamics. 
We assume that the nonlinear term, $f(\cdot,\cdot,t)$, in both agents' and reference inputs' dynamics satisfies the Lipschitz-type condition.	
First, a local filter is introduced for each agent to estimate the average of the reference inputs and the reference velocities.
Then, the control input $u_i$, $i=1,\ldots,n$, is designed for each agent such that $x_i$ and $v_i$ track, respectively, $p_i$ and $q_i$, where $p_i$, $q_i \in \mathbb{R}^\mathfrak{p}$ are the filter outputs.
	For notational simplicity, we will remove the index $t$ from variables in the reminder of the paper.\newline	
\begin{assumption}\label{lip} \cite{yu2010second}
			The vector-valued term $f(\cdot,\cdot,t)$ satisfies the following inequality $\forall t \geq 0$
			\begin{align*}
			\|f(x,v,t)-f(y,z,t)\|_1 \leq  \rho_1 \|x-y\|_1  +\rho_2 \|v-z\|_1,  
			\end{align*}
			where $x$, $v$, $y$, $z \in \mathbb{R}^\mathfrak{p}$, $\rho_1$, $\rho_2 \in \mathbb{R}^+$ and $f(0,0,t)=0$.
		\end{assumption}
\begin{remark}
Note that Assumption \ref{lip} is a Lipschitz-type condition, satisfied by many well-known systems such as the pendulum system with a control torque, car-like robots, the Chua’s circuit, the Lorenz system, and the Chen system \cite{mei2013distributed}.
\end{remark}
Under Assumption \ref{lip}, the unknown term $f(\cdot,\cdot,t)$ might be unbounded. Therefore, two novel state-dependent time varying gains will be introduced to overcome the unboundedness effect of this term.


The following local filter is introduced for agent $i$, ${i=1,\ldots,n}$, to estimate the average of reference inputs and reference velocities 
\begin{align}\label{p-term2}
	\small
	p_i=& z_i+r_i, \notag \\
	\ddot{z}_i=&- \kappa (p_i-r_i)- \kappa (q_i-v^r_i) \notag \\
	&  -\alpha \psi_i   \mbox{sgn} \Big[  \sum\limits_{j=1}^{n} a_{ij} \Big \{ (p_i+q_i)-(p_j+q_j) \Big \} \Big],
\end{align}
\normalsize
where $q_i=\dot{p}_i$,  $z_i \in \mathbb{R}^\mathfrak{p}$ is an auxiliary filter variable, $\psi_i=\|r_i\|_1+\|v^r_i\|_1+\gamma$, and $\kappa$, $\alpha$, $\gamma \in \mathbb{R}^+$ are control gains to be designed.
The distributed control input $u_i$, $i=1,\ldots,n$, is designed to drive $x_i$ and $v_i$ to track, respectively, $p_i$ and $q_i$,
\begin{align}\label{con-inp2}
	u_i=&-\eta \tilde{x}_i-\eta \tilde{v}_i \notag -\eta (\|x_i-r_i\|_1+\|v_i-v^r_i\|_1+\gamma) \times \notag \\ &\mbox{sgn}(\tilde{x}_i+\tilde{v}_i)  +\ddot{z}_i,  
\end{align}
where $\tilde{x}_i=x_i-p_i$, $\tilde{v}_i=v_i-q_i$ and $\eta \in \mathbb{R}^+$ is a control gain to be designed.
\begin{theorem} \label{thm:DAC-lip}
	Under the control law given by  \eqref{p-term2} and \eqref{con-inp2} for \eqref{non-agents-dynamic}, the distributed average tracking goal \eqref{goal} is achieved asymptotically, provided that Assumptions \ref{conn-graph} and \ref{lip} hold and the control gains $\alpha$, $\gamma$, and $\eta$ are chosen such that $\kappa>1$, $\alpha >  \max \{\rho_1,\rho_2 \}+\kappa$ and $\eta > \max \{ 1, \rho_1,\rho_2 \}$, where $\rho_1$ and $\rho_2$ are defined in Assumption \ref{lip}.
\end{theorem}

\emph{Proof}:
The proof contains two steps. 
First, it is proved that for the $i$th agent, $\lim\limits_{t \to \infty} p_i = \frac{1}{n} \sum\limits_{j=1}^{n} r_j$ and $\lim\limits_{t \to \infty} q_i = \frac{1}{n} \sum\limits_{j=1}^{n} v^r_j$.
Second, it is shown that by using the control input \eqref{con-inp2} for agent $i$, $\lim\limits_{t \to \infty} x_i = p_i$ and $\lim\limits_{t \to \infty} v_i = q_i$.

Using $q_i=\dot{p}_i$, the filter dynamics \eqref{p-term2} can be rewritten as
\begin{align}\label{clos-lip}
\small
	\dot{p}_i=& q_i, \notag \\
	\dot{q}_i=&- \kappa (p_i-r_i)- \kappa (q_i-v^r_i) \notag \\
	&  -\alpha \psi_i   \mbox{sgn} \Big[  \sum\limits_{j=1}^{n} a_{ij} \{ (p_i+q_i)-(p_j+q_j) \} \Big] \notag \\
	& + f(r_i,v^r_i,t),
\end{align}
\normalsize
Let $\tilde{p}=(M \otimes I_\mathfrak{p}) p$ and $\tilde{q}=(M \otimes I_\mathfrak{p}) q$, where $p=[p_1^T,\ldots,p_n^T]^T$, $q=[q_1^T,\ldots,q_n^T]^T$ and $M=I_n- \frac{1}{n}\textbf{1}^T_n \textbf{1}_n$.
Now the local filter's closed dynamics \eqref{clos-lip} can be rewritten in vector form as
\begin{align*}
\small
	\dot{\tilde{p}}=& \tilde{q}, \notag \\
	\dot{\tilde{q}}=& -\kappa \tilde{p} +\kappa  (M  \otimes I_\mathfrak{p}) r - \kappa \tilde{q} +\kappa (M  \otimes I_\mathfrak{p}) v^r \notag \\
	& -\alpha (M  \psi \otimes I_\mathfrak{p}) \mbox{sgn} [( L \otimes I_\mathfrak{p})(\tilde{p}+\tilde{q})] \notag \\
	&+(M \otimes I_\mathfrak{p})  f(r,v^r,t),
\end{align*}
\normalsize
where $\psi=\mbox{diag}(\psi_1,\ldots,\psi_n)$. 
Consider the following Lyapunov function candidate
$V_1= \frac{1}{2} 
	\begin{bmatrix}
		\tilde{p}^T && \tilde{q}^T
	\end{bmatrix}
	( L \otimes
	\begin{bmatrix}
		2 \kappa  && 1 \\
		1  && 1
	\end{bmatrix}
	\otimes I_\mathfrak{p})
	\begin{bmatrix}
		\tilde{p} \\
		\tilde{q}
	\end{bmatrix}$.
Since $(\mathbf{1}_n \otimes I_\mathfrak{p}) ^T \tilde{p} =\mathbf{0}_{n\mathfrak{p}}$ and $(\mathbf{1}_n \otimes I_\mathfrak{p}) ^T \tilde{q} =\mathbf{0}_{n\mathfrak{p}}$, by using Lemma \ref{eigen}, we have
$V_1 \geq \frac{\lambda_2(L)}{2} \begin{bmatrix}
		\tilde{p}^T && \tilde{q}^T
	\end{bmatrix} (
	\begin{bmatrix}
		2 \kappa && 1 \\
		1 && 1
	\end{bmatrix} \otimes I_{n\mathfrak{p}})
	\begin{bmatrix}
		\tilde{p} \\
		\tilde{q}
	\end{bmatrix}$,
where $\lambda_2 (L)$ is defined in Lemma \ref{eigen}.
It can be proved that if $\kappa >\frac{1}{2}$, then $\begin{bmatrix}
		2 \kappa && 1 \\
		1 && 1
\end{bmatrix}>0$, which means that $V_1$ is positive definite.
The derivative of $V_1$ is given as
\begin{align*}
	\small
	\dot{V}_1 =&  2 \kappa \tilde{p}^T (L \otimes I_\mathfrak{p}) \tilde{q}+\tilde{q}^T (L \otimes I_\mathfrak{p}) \tilde{q}-\kappa \tilde{p}^T (L \otimes I_\mathfrak{p}) \tilde{p} \notag \\
	&+\kappa \tilde{p}^T (L \otimes I_\mathfrak{p})r -\kappa \tilde{p}^T (L \otimes I_\mathfrak{p}) \tilde{q} +\kappa \tilde{p}^T (L \otimes I_\mathfrak{p}) v^r \notag \\
	&-\alpha\tilde{p}^T (L  \psi \otimes I_\mathfrak{p}) \mbox{sgn} [( L \otimes I_\mathfrak{p})(\tilde{p}+\tilde{q})] \notag \\
	&+\tilde{p}^T (L \otimes I_\mathfrak{p})  f(r,v^r,t)-\kappa \tilde{q}^T (L \otimes I_\mathfrak{p}) \tilde{p} \notag \\
	&+\kappa \tilde{q}^T (L \otimes I_\mathfrak{p})r -\kappa \tilde{q}^T (L \otimes I_\mathfrak{p}) \tilde{q} +\kappa \tilde{q}^T (L \otimes I_\mathfrak{p}) v^r \notag \\
	&-\alpha\tilde{q}^T (L  \psi \otimes I_\mathfrak{p}) \mbox{sgn} [( L \otimes I_\mathfrak{p})(\tilde{p}+\tilde{q})] \notag \\
	&+\tilde{q}^T (L \otimes I_\mathfrak{p})  f(r,v^r,t)\\
	=&-\kappa \tilde{p}^T (L \otimes I_\mathfrak{p}) \tilde{p} - (\kappa -1 ) \tilde{q}^T (L \otimes I_\mathfrak{p}) \tilde{q} \\
	&+ \kappa (\tilde{p}+\tilde{q})^T (L \otimes I_\mathfrak{p}) (r + v^r) \\
	& -\alpha \sum\limits_{i=1}^{n} \psi_i \Big [ \sum\limits_{j=1}^{n} a_{ij} \Big \{ (\tilde{p}_i+\tilde{q}_i)-(\tilde{p}_j+\tilde{q}_j) \Big\} \Big ]^T \times \\
	& \mbox{sgn} \Big[ \sum\limits_{j=1}^{n} a_{ij} \Big\{ (\tilde{p}_i+\tilde{q}_i)-(\tilde{p}_j+\tilde{q}_j) \Big\} \Big] \\
	& + \sum\limits_{i=1}^{n} \Big [ \sum\limits_{j=1}^{n} a_{ij} \Big\{(\tilde{p}_i+\tilde{q}_i)-(\tilde{p}_j+\tilde{q}_j) \Big\} \Big ]^T f(r_i,v^r_i,t) \\
	=&-\kappa \tilde{p}^T (L \otimes I_\mathfrak{p}) \tilde{p} - (\kappa -1 ) \tilde{q}^T (L \otimes I_\mathfrak{p}) \tilde{q} \\
	&+ \kappa \sum\limits_{i=1}^{n} \Big [ \sum\limits_{j=1}^{n} a_{ij} \Big\{(\tilde{p}_i+\tilde{q}_i)-(\tilde{p}_j+\tilde{q}_j) \Big \} \Big ]^T \Big[ r_i+v^r_i \Big] \\
	& -\alpha \sum\limits_{i=1}^{n} \psi_i \Big [ \sum\limits_{j=1}^{n} a_{ij} \Big \{ (\tilde{p}_i+\tilde{q}_i)-(\tilde{p}_j+\tilde{q}_j) \Big\} \Big ]^T \times \\
	& \mbox{sgn} \Big[ \sum\limits_{j=1}^{n} a_{ij} \Big\{ (\tilde{p}_i+\tilde{q}_i)-(\tilde{p}_j+\tilde{q}_j) \Big\} \Big] \\
	&+\sum\limits_{i=1}^{n} \Big [ \sum\limits_{j=1}^{n} a_{ij} \Big \{(\tilde{p}_i+\tilde{q}_i)-(\tilde{p}_j+\tilde{q}_j) \Big \} \Big ]^T \times \\
	& \Big[ f(r_i,v^r_i,t)-f(0,0,t) \Big],
\end{align*}
\normalsize
where $\tilde{p}_i$ and $\tilde{q}_i$ are, respectively, the $i$th components of $\tilde{p}$ and $\tilde{q}$ and we have used $LM=L$ and Assumption \ref{lip} to obtain, respectively, the first and the third equalities.
Under Assumption \ref{lip} and using the triangular inequality, we have
\begin{align*}
\small
	\dot{V}_1 \leq&-\kappa \tilde{p}^T (L \otimes I_\mathfrak{p}) \tilde{p} - (\kappa -1 ) \tilde{q}^T (L \otimes I_\mathfrak{p}) \tilde{q} \\
	&+\kappa \sum\limits_{i=1}^{n} \Big \| \sum\limits_{j=1}^{n} a_{ij} \Big \{ (\tilde{p}_i+\tilde{q}_i)-(\tilde{p}_j+\tilde{q}_j) \Big \} \Big \|_1 \times \\
	& ( \|r_i\|_1+\|v^r_i\|_1) \\
	&-\alpha \sum\limits_{i=1}^{n} \psi_i \Big \|   \sum\limits_{j=1}^{n} a_{ij} \Big \{ (\tilde{p}_i+\tilde{q}_i)-(\tilde{p}_j+\tilde{q}_j) \Big \} \Big \|_1 \\
	&+\sum\limits_{i=1}^{n} \Big \| \sum\limits_{j=1}^{n} a_{ij} \Big \{ (\tilde{p}_i+\tilde{q}_i)-(\tilde{p}_j+\tilde{q}_j) \Big \}  \Big \|_1 \times \\
	&  (\rho_1 \|r_i\|_1+\rho_2\|v^r_i\|_1) \\
	= & -\kappa \tilde{p}^T (L \otimes I_\mathfrak{p}) \tilde{p} - (\kappa -1 ) \tilde{q}^T (L \otimes I_\mathfrak{p}) \tilde{q} \\
	&-\alpha \sum\limits_{i=1}^{n} \psi_i \Big \|   \sum\limits_{j=1}^{n} a_{ij} \Big \{ (\tilde{p}_i+\tilde{q}_i)-(\tilde{p}_j+\tilde{q}_j) \Big \} \Big \|_1 \\
	&+\sum\limits_{i=1}^{n}  \Big ( (\kappa+\rho_1) \|r_i\|_1+(\kappa+\rho_2)\|v^r_i\|_1 \Big )  \times \\
	& \Big \| \sum\limits_{j=1}^{n} a_{ij} \Big \{ (\tilde{p}_i+\tilde{q}_i)-(\tilde{p}_j+\tilde{q}_j) \Big \}  \Big \|_1 \\
	= & -\kappa \tilde{p}^T (L \otimes I_\mathfrak{p}) \tilde{p} - (\kappa -1 ) \tilde{q}^T (L \otimes I_\mathfrak{p}) \tilde{q} \\
	&+\sum\limits_{i=1}^{n} \Big ( (\kappa+\rho_1-\alpha) \|r_i\|_1+(\kappa+\rho_2-\alpha)\|v^r_i\|_1 -\alpha \gamma\Big )  \times \\
	& \Big \| \sum\limits_{j=1}^{n} a_{ij} \Big \{ (\tilde{p}_i+\tilde{q}_i)-(\tilde{p}_j+\tilde{q}_j) \Big \}  \Big \|_1 ,
\end{align*}
\normalsize
where we have used the definition of $\psi_i$ to obtain the last equality.
Since $\alpha >  \max \{\rho_1,\rho_2 \}+\kappa$, we will have
\begin{align*}
\dot{V}_1 \leq & -\kappa \tilde{p}^T (L \otimes I_\mathfrak{p}) \tilde{p} - (\kappa -1 ) \tilde{q}^T (L \otimes I_\mathfrak{p}) \tilde{q} \\ 
\leq & -\kappa \lambda_2(L) \tilde{p}^T  \tilde{p}-(\kappa -1)\lambda_2(L) \tilde{q}^T  \tilde{q} < 0,
\end{align*}
where we have used Lemma \ref{eigen}, since $(\mathbf{1}_n \otimes I_\mathfrak{p}) ^T \tilde{p} =\mathbf{0}_{n\mathfrak{p}}$ and $(\mathbf{1}_n \otimes I_\mathfrak{p}) ^T \tilde{q}=\mathbf{0}_{n\mathfrak{p}}$, and  $\kappa>1$ to obtain the second inequality.
Using Theorem 4.10 in \cite{khalil2002nonlinear}, it is concluded that $\begin{bmatrix}
\tilde{p} \\
\tilde{q}
\end{bmatrix}=\mathbf{0}_{2n \mathfrak{p}}$ is globally exponentially stable, which means $ \lim\limits_{t \to \infty} p_i = \frac{1}{n} \sum\limits_{j=1}^{n} p_j$ and $ \lim\limits_{t \to \infty} q_i = \frac{1}{n} \sum\limits_{j=1}^{n} q_j$ for $i=1,\ldots,n$.
Now it is enough to show that $\lim\limits_{t \to \infty} \sum\limits_{j=1}^{n} p_j=  \sum\limits_{j=1}^{n} r_j$ and $ \lim\limits_{t \to \infty} \sum\limits_{j=1}^{n} q_j=\sum\limits_{j=1}^{n} v^r_j$ which results in $\lim\limits_{t \to \infty } p_i = \frac{1}{n} \sum\limits_{j=1}^{n} r_j$ and $\lim\limits_{t \to \infty } q_i = \frac{1}{n} \sum\limits_{j=1}^{n} v^r_j$.
Define the variables ${S_1=\sum_{i=1}^n (p_i- r_i)}$ and ${S_2=\sum_{i=1}^n (q_i-v_i^r) }$, we can get from \eqref{clos-lip} that
\begin{align}\label{sum-lip}
\dot{S}_1=& S_2, \notag \\
\dot{S}_2=&- \kappa S_1- \kappa S_2 \notag \\
&  -\alpha \sum\limits_{i=1}^{n} \psi_i   \mbox{sgn} \Big[ \sum\limits_{j=1}^{n} a_{ij} \Big\{ (p_i+q_i)-(p_j+q_j) \Big\} \Big],
\end{align}
We then use input-to-state stability to analyze the system \eqref{sum-lip} by treating the term $\sum\limits_{i=1}^{n} \psi_i   \mbox{sgn} \Big[ \sum\limits_{j=1}^{n} a_{ij} \Big\{ (p_i+q_i)-(p_j+q_j) \Big\} \Big]$ as the input and $S_1$ and $S_2$ as the states.
Since $\kappa>1$, the matrix $\begin{bmatrix}
\mathbf{0}_\mathfrak{p} & I_\mathfrak{p} \\
-\kappa I_\mathfrak{p} & -\kappa I_\mathfrak{p}
\end{bmatrix}$ is Hurwitz. 
Thus, the system \eqref{sum-lip} with zero input is exponentially stable and hence input-to-state stable. 
Since $p_i+q_i \to p_j+q_j$, $i,j=1,\cdots,n$, as $t \to \infty$, it follows that $S_1 \to 0$ and $S_2 \to 0$.
Therefore, we have that $\lim\limits_{t \to \infty}\sum_{i=1}^n p_i = \sum_{i=1}^n r_i$ and $\lim\limits_{t \to \infty} \sum_{i=1}^n q_i = \sum_{i=1}^n v_i^r$, respectively. 
Employing the results of these two parts, it is concluded that $\lim\limits_{t \to \infty} p_i \to \frac{1}{n} \sum_{j=1}^n r_j$ and ${\lim\limits_{t \to \infty} q_i \to \frac{1}{n} \sum_{j=1}^n v_j^r}$.

In the remaining, we will prove that $\lim\limits_{t \to \infty} x_i = p_i$ and $\lim\limits_{t \to \infty} v_i = q_i$ asymptotically.
Using the control input \eqref{con-inp2} for \eqref{non-agents-dynamic}, we get the closed-loop dynamics in vector form as
\begin{align*}
\small
	\dot{\tilde{x}}=&\tilde{v}, \\
	\dot{\tilde{v}}=& f(x,v,t) -\eta \tilde{x}-\eta \tilde{v} \\
	&-\eta (\|x-r\|_1+\|v-v^r\|_1+\gamma) \mbox{sgn}(\tilde{x}+\tilde{v}) - f(r,v^r,t), 
\end{align*}
\normalsize
where $\tilde{x}=[\tilde{x}_1^T,\ldots,\tilde{x}_n^T]^T$, $\tilde{v}=[\tilde{v}_1^T,\ldots,\tilde{v}_n^T]^T$ and $f(x,v,t)=[f^T(x_1,v_1,t),\ldots,f^T(x_n,v_n,t)]^T$.
Consider the candidate Lyapunov function
$V_2= \frac{1}{2} 
	\begin{bmatrix}
		\tilde{x}^T && \tilde{v}^T
	\end{bmatrix} 
	\begin{bmatrix}
		2\eta I_{n\mathfrak{p}} && I_{n\mathfrak{p}} \\
		I_{n\mathfrak{p}} && I_{n\mathfrak{p}}
	\end{bmatrix}
	\begin{bmatrix}
		\tilde{x} \\
		\tilde{v}
	\end{bmatrix}$.
Since $\eta>\frac{1}{2}$, $V_2$ is positive definite. By taking the derivative of $V_2$, we will have
\begin{align*}
	\dot{V}_2=&  2\eta \tilde{x}^T  \tilde{v}+\tilde{v}^T \tilde{v} -\eta \tilde{x}^T (\tilde{x}+\tilde{v})\\
	& -\eta \sum\limits_{i=1}^{n} (\|x_i-r_i\|_1+\|v_i-v^r_i\|_1+\gamma)\tilde{x}_i^T\mbox{sgn}(\tilde{x}_i+\tilde{v}_i) \\ &+\sum\limits_{i=1}^{n} \tilde{x}_i^T (f(x_i,v_i,t)-f(r_i,v^r_i,t))
	-\eta \tilde{v}^T (\tilde{x}+\tilde{v}) \\
	& -\eta \sum\limits_{i=1}^{n} (\|x_i-r_i\|_1+\|v_i-v^r_i\|_1+\gamma) \tilde{v}_i^T\mbox{sgn}(\tilde{x}_i+\tilde{v}_i) \\ 
	&+\sum\limits_{i=1}^{n} \tilde{v}_i^T (f(x_i,v_i,t)-f(r_i,v^r_i,t)) \\
	=&-\eta \tilde{x}^T  \tilde{x}+(1-\eta) \tilde{v}^T  \tilde{v} \\
	&-\eta \sum\limits_{i=1}^{n} (\|x_i-r_i\|_1+\|v_i-v^r_i\|_1+\gamma)\|\tilde{x}_i+\tilde{v}_i\|_1 \\ &+\sum\limits_{i=1}^{n}(\tilde{x}_i+\tilde{v}_i)^T (f(x_i,v_i,t)-f(r_i,v^r_i,t)) \\
	\leq & -\eta \tilde{x}^T  \tilde{x}+(1-\eta)\tilde{v}^T  \tilde{v}-\eta \gamma \sum\limits_{i=1}^{n} \|\tilde{x}_i+\tilde{v}_i\|_1  \\
	&-\eta \sum\limits_{i=1}^{n} (\|x_i-r_i\|_1+\|v_i-v^r_i\|_1)\|\tilde{x}_i+\tilde{v}_i\|_1 \\
	& +\sum\limits_{i=1}^{n} (\rho_1 \|x_i-r_i\|_1+\rho_2 \|v_i-v^r_i\|_1 ) \|\tilde{x}_i+\tilde{v}_i\|_1,
\end{align*}
where we have used Assumption \ref{lip} to obtain the inequality.
Since $\eta >\max \{  \rho_1,\rho_2 \}$, we can get that 
$\dot{V}_2 	\leq  -\eta \tilde{x}^T  \tilde{x} -(\eta-1) \tilde{v}^T  \tilde{v} \leq 0$,
where we have used $\eta >1$ to obtain the last inequality.
Using Theorem 4.10 in \cite{khalil2002nonlinear}, it is concluded that $\begin{bmatrix}
\tilde{x} \\
\tilde{v}
\end{bmatrix}=\mathbf{0}_{2n\mathfrak{p}}$ is globally exponentially stable. Therefore, it is concluded that for $i$th agent, $i=1,\ldots,n$, $\lim\limits_{t \to \infty} x_i = \frac{1}{n} \sum\limits_{j=1}^{n} r_j$ and 
$ \lim\limits_{t \to \infty} v_i = \frac{1}{n} \sum\limits_{j=1}^{n} v^r_j$ asymptotically.
\endproof
\begin{remark} 
	Due to the presence of the unknown term $f(\cdot,\cdot,t)$ in the agents' dynamics, the proposed algorithms for double-integrator agents are not applicable to achieve the distributed average tracking. 
	For example, by employing the algorithm in \cite{feidoubleintegrator} for \eqref{non-agents-dynamic}, the two equalities $\sum_{j=1}^n x_j =\sum_{j=1}^n r_j$ and $\sum_{j=1}^n v_j = \sum_{j=1}^n v_j^r$ do not hold anymore.
	In fact, the unknown term $f(\cdot,\cdot,t)$ functions as a disturbance and it will not allow the average of the positions and velocities to track the reference inputs and the reference velocities, respectively.
	This shows the essence of using the local filter \eqref{z-term} in our algorithm.
\end{remark}
\subsection{Discussion}\label{section1}
	Suppose that the unknown term $f(\cdot,\cdot,t)$ satisfies the boundedness condition
	\begin{align}\label{bounded-f}
	{\text{sup}_{t \in [0,\infty)} \| f(x,v,t)\|_1 \leq \bar{f}},
	\end{align}
	 where $x$, $v \in \mathbb{R}^\mathfrak{p}$, and $\bar{f} \in \mathbb{R}^+$.
	Then, the proposed algorithm in Subsection \ref{section2} can be simplified as the following local filter dynamics and control input with constant control gains to achieve the distributed average tracking.
	\begin{align}
		p_i=& z_i+r_i, 
		\notag	\\
		\ddot{z}_i=&-\alpha  \sum\limits_{j=1}^{n} a_{ij} \mbox{sgn}\Big[(p_i+q_i)-(p_j+q_j) \Big], \label{z-term} \\
	u_i=& -\eta \mbox{sgn} \Big [\tilde{x}_i+\tilde{v}_i \Big] +\ddot{z}_i,
	\label{con-inp}
	\end{align}	
	where $p_i$, $q_i$, $z_i$, $\tilde{x}_i$ and $\tilde{v}_i$ are defined in Subsection \ref{section2} and  $\alpha$, $\eta \in \mathbb{R}^+$ are control gains to be defined. 
	\begin{assumption}\label{initial}
		The variables $z_i$ and $\dot{z}_i$ are initialized such that $\sum\limits_{i=1}^{n} z_i(0)=\mathbf{0}_{\mathfrak{p}}$, and $\sum\limits_{i=1}^{n} \dot{z}_i(0)=\mathbf{0}_{\mathfrak{p}}$.	\footnote{A special choice is $z_i(0)=\mathbf{0}_{\mathfrak{p}}$ and $\dot{z}_i(0)=\mathbf{0}_{\mathfrak{p}}$ for all ${i=1,\ldots,n}$.} 
	\end{assumption}
	\begin{theorem} \label{DAC-non-dynamic}
		Under the control law given by \eqref{z-term} and \eqref{con-inp} for \eqref{non-agents-dynamic}, the distributed average tracking goal \eqref{goal} is achieved asymptotically, provided that Assumption \ref{conn-graph} and the boundedness condition \eqref{bounded-f} hold and the control gains $\alpha$ and $\eta$ are chosen such that
		$\alpha >  \frac{2 \|s(0)\|_1 + \|(D^T \otimes I_\mathfrak{p}) x(0)\|_1 + n\bar{f}}{\lambda_2 (L)}$ and $\eta >  2\|\tilde{s}(0)\|_1 + \|\tilde{x}(0)\|_1 +2 \bar{f}$, where $\lambda_2 (L)$ and $\bar{f}$ are defined, respectively, in Lemma \ref{eigen} and  \eqref{bounded-f}.
	\end{theorem}
	\emph{Proof}:
	Similar to the Theorem \ref{thm:DAC-lip}, the proof contains two steps.  
	The general idea of the first step of the proof is adopted from  \cite{feidoubleintegrator}, where $a_i^r$ is replaced by $f(r_i,v^r_i,t)$.
	The second step of proof is similar to the second step of proof of Theorem \ref{thm:DAC-lip}.
	Thus, the proof is omitted 
	\endproof
\begin{remark} 
	In \eqref{z-term}-\eqref{con-inp}, each agent needs its neighbors' filters' outputs besides its own filter outputs, absolute position and absolute velocity, while in \eqref{p-term2}-\eqref{con-inp2}, the agent needs its neighbors' reference inputs and reference velocities too.
	Further, in both algorithms, there is no need for correct position and velocity initialization, where the initialization of the physical variables might not be feasible for real applications. 
	In the algorithm \eqref{z-term}-\eqref{con-inp}, only the initialization of the filters’ auxiliary variables is required, which can be easily satisfied.
\end{remark}

\begin{remark}
	The unknown bounded term $f(\cdot,\cdot,t)$ that satisfies \eqref{bounded-f}, can be interpreted as a perturbation resulting from modeling errors, uncertainties, and disturbances that exist in many realistic problems.
	However, Assumption \ref{lip} includes a more general group of systems including unbounded nonlinear terms. 
	In addition if $\gamma$ is chosen properly, the algorithm \eqref{p-term2}-\eqref{con-inp2} is still applicable in the presence of an additive bounded term in either agents' or reference inputs' dynamics.
\end{remark}
\section{CONCLUSIONS}
In this paper, the distributed average tracking of physical second-order agents with nonlinear dynamics was studied. 
A nonlinear term was assumed in both agents' and reference inputs' dynamics satisfying the Lipschitz-type condition. 
Due to the presence of the nonlinear term in agents' dynamics, a filter design combined with controller design was introduced to solve the problem. 
The idea was that each filter outputs converge to the average of the reference inputs and the reference velocities asymptotically and in parallel the agent's position and velocity are driven to track the filter outputs. 
Since the nonlinear term could be unbounded, novel state-dependent time varying gains were employed in each agent's filter and control input to overcome the unboundedness effect.

\bibliographystyle{IEEEtran}
\bibliography{refs}

\end{document}